\documentclass[11pt]{article}
\advance \textheight by -1 \topmargin
\advance \textheight by -1 \headheight
\advance \textheight by -1 \headsep
\advance \textheight by -2 \footskip
\vsize = \textheight

\hsize = \textwidth
\usepackage{amssymb,amsmath,courier,array}

\newtheorem{theorem}{Theorem}[section] 

\newtheorem{lemma}[theorem]{Lemma}
\newtheorem{proposition}[theorem]{Proposition}

\newtheorem{remark}[theorem]{Remark}

\begin{document}

\newcommand{\proof}{{\bf Proof.~}}
\newcommand{\qed}{\Box}

\renewcommand{\bar}{\overline}
\renewcommand{\leq}{\leqslant}
\renewcommand{\geq}{\geqslant}
\newcommand{\la}{\langle}
\newcommand{\ra}{\rangle}
\newcommand{\aut}{{\rm Aut}}
\newcommand{\inn}{{\rm Inn}}
\newcommand{\out}{{\rm Out}}
\newcommand{\Z}{\Bbb{Z}}
\newcommand{\qq}{{\bf q}} 
\newcommand{\bb}{{\bf b}} 
\newcommand{\QQ}{{\bf Q}} 
\newcommand{\BB}{{\bf B}}


\title{On groups with a class-preserving outer automorphism}
\author{Peter A. Brooksbank~~~~~~~~Matthew S. Mizuhara}

\maketitle

\begin{abstract}
In 1911, Burnside asked whether or not there exist groups that have an outer
automorphism which preserves conjugacy classes. Two years later he answered
his own question by constructing a family of such groups. Using the small group
library in {\sc magma} we determine all of the groups of order $n<512$ 
that possess such an automorphism. Our investigations led to the discovery of four new 
infinite families of such groups, all of which are 2-groups of coclass 4.
\end{abstract}


\section{Introduction}
\label{intro}
Let $G$ be a group, $\aut(G)$ the automorphism group of $G$, 
and $\inn(G)$ the subgroup of inner automorphisms.
Then $\aut(G)$ acts naturally on the set of conjugacy classes of $G$, and we 
denote the kernel of this action by $\aut_c(G)$. Evidently $\inn(G)\leq\aut_c(G)$,
and we accordingly refer to the elements of $\aut_c(G)$ as 
{\em nearly inner automorphisms} or {\em class-preserving automorphisms}.

In 1911, Burnside~\cite[Note B]{burnside1} asked the question: {\em Are there groups 
$G$ having nearly inner automorphisms that are not inner?} In 1913, he answered his 
own question~\cite{burnside2}: for each prime $p\equiv\pm 3\;({\rm mod}\;8)$, 
there is such a group of order 
$p^6$ and class 2.

With the aid of modern computer technology we sought to investigate how rare this
property is among groups of small order. Searching the small group
library (made available by the work of Besche, Eick and O'Brien~\cite{millenium})
in the {\sc magma} system~\cite{magma} 
we determined, for each $n<512$, 
all groups $G$ of order $n$ such that 
$\out_c(G):=\aut_c(G)/\inn(G)\neq 1$. The algorithmic techniques we used to carry out 
this search, and some of the data that we generated from it, are presented in the concluding 
section of the paper. 

One striking feature of the data is that the groups Burnside sought appear to be far more 
common than one might have thought. For example, over $60\%$ of the 56,092 groups of 
order 256 have a nearly inner automorphism that is not inner. (General theory suggests that 
one should look for such groups among soluble groups; Feit and 
Seitz~\cite[Section C]{feit-seitz} showed that if $G$ is a finite simple group then 
$\out_c(G)=1$.)

Several suggestive patterns emerged during our systematic search for groups $G$ with $\out_c(G)\neq 1$.
Closer examination led to the discovery of new infinite families of such groups. Indeed the main result
that we report in this paper is the following.


\begin{theorem}
\label{thm:main}
There are four distinct infinite families 
$\displaystyle{\mathcal{H}=\{H_{j}\}_{j=1}^{\infty}}$,
where $H_{j}$ is a 4-generator 2-group of order
$2^{5+j}$ and class $j+1$ such that $\out_c(H_{j})\neq 1$.
\end{theorem}

Since Burnside's initial discovery, the problem has been revisited on many occasions. 
We briefly outline results that relate to Theorem~\ref{thm:main}, 
and refer the reader to~\cite{yadav} for a more comprehensive survey of the problem
and its applications.
\smallskip

In 1947, Wall showed that, for each integer $m$ divisible by 8, the general linear group 
${\rm GL}(1,\Z/m)$, of order $m\cdot\varphi(m)$ has a nearly inner automorphism that is 
not inner~\cite{wall}. These include a smallest example of such groups, namely 
${\rm GL}(1,\Z/8)$, of order 32.  (There are actually two non-isomorphic groups of order 
32 having this property.) The 2-groups in Wall's family, namely ${\rm GL}(1,\Z/2^k)$,  
have class 2.

In 1979, Heineken constructed, for each odd prime $p$, an infinite family of 
$p$-groups of class 2, {\em all} of whose automorphisms are 
class-preserving~\cite{heineken}. As far as we are aware, these are the only 
known infinite families of groups $G$ for which $\aut_c(G)=\aut(G)$.
Our search revealed that the smallest such group has order $2^7$; the 
smallest group in Heineken's family has order $3^{15}$.

In 2001, Hertweck constructed a family of Frobenius groups as subgroups
of affine semi-linear groups $A\Gamma(F)$, where $F$ is a finite field, which 
possess nearly inner automorphisms that are not inner~\cite{hertweck}. 
The smallest group in Hertweck's family has order 300.

Until fairly recently, all of the known infinite families of groups having
class-preserving outer automorphisms were nilpotent groups of class 2.
In 1992, Malinowska~\cite{malinowska} exhibited, for each prime $p>5$ 
and each $r>2$, a $p$-group $G$ of class $r$ such that $\out_c(G)\neq 1$;
it is not clear how the order of $G$ grows as a function of $r$.
\smallskip

It is evident from the statement of Theorem~\ref{thm:main} that the order and
class of the groups $H_j$ in each family grow in a very controlled way.
This is because $H_{j+1}$ is built as a central extension of $H_j$ by $\Z/2$.
Indeed, each $\mathcal{H}$ may be constructed algorithmically using
the $p$-group generation algorithm~\cite{obrien}; this is 
precisely how the families were discovered and studied. Furthermore, the groups
in all four families have coclass 4, so we have shown that they are all
``mainline groups" in the coclass graph $\mathcal{G}(2,4)$ (cf.~\cite{eick-LG}).
\smallskip

Briefly, the paper is organized as follows. In Section~\ref{prelims} we summarize
the necessary background on $p$-groups. The families alluded to in 
Theorem~\ref{thm:main} are introduced in Section~\ref{four-families}. These are
parametrized by vectors $\epsilon\in\{0,1\}^4$, but each $\mathcal{H}^{\epsilon}$
is isomorphic to three others, yielding four distinct families. The proof of 
Theorem~\ref{thm:main} is given in Section~\ref{proof}. Finally, in Section~\ref{data}, 
we describe the elementary algorithm that we used to search the small group database, 
and summarize the data obtained from this search.

\section{Preliminaries}
\label{prelims}
In this short section we introduce some standard notation and terminology,
and summarize the necessary theory. 

If $x,y$ are elements of a group
$G$, then we write $x^y=y^{-1}xy$ and $[x,y]=x^{-1}x^y$.
If $X,Y$ are subsets of $G$, then $[X,Y]=\la [x,y]\colon x\in X,y\in Y\ra$.
If $H=\la X\ra$ and $K=\la Y\ra$, then $[H,K]=[X,Y]^{\la H,K\ra}$. The
{\em lower central series} of $G$ is the series
\begin{equation}
\label{eq:lower}
G=\gamma_1(G)\geq\gamma_2(G)\geq\ldots
\end{equation}
where $\gamma_{i+1}(G)=[G,\gamma_i(G)]$. A group $G$ is {\em nilpotent}
if $\gamma_j(G)=1$ for some $j\geq 1$, in which case the smallest $c$ such  
that $\gamma_{c+1}(G)=1$
is called the {\em nilpotency class} (or simply {\em class}) of $G$.
A finite group $G$ is a {\em $p$-group} if $|G|=p^n$ for some prime $p$.
All $p$-groups are nilpotent, and if $G$ has class $c$, then $G$ has
{\em coclass} $n-c$. If a $p$-group $G$ is minimally generated by $d$
elements, then we say that $G$ is a {\em $d$-generator} group.
\medskip

Each nilpotent group (more generally, each soluble group) possesses 
a {\em polycyclic generating sequence}~\cite[Chapter 8]{handbook}.
This in turn gives rise to a {\em power-conjugate presentation} (or simply
{\em pc-presentation}), an extremely efficient model for
computing with soluble groups. We describe
these presentations specifically for $p$-groups.
\smallskip

Let $G$ be a $p$-group and 
$X=[x_1,\ldots,x_n]\subset G$ be 
such that if $P_i=\la x_i,\ldots,x_n\ra$ ($i=1,\ldots,n$), then 
$P_i/P_{i+1}$ has order $p$, and $G=P_1>P_2>\ldots>P_n>1$
refines the lower central series in Equation~(\ref{eq:lower}). 
If $G$ has nilpotency class $c$, we define a {\em weighting},
$w\colon X\to\{1,\ldots,c\}$, where $w(x_i)=k$
if $x_i\in \gamma_{k-1}(G)\setminus\gamma_k(G)$. Evidently,
$w(x_i)\geq w(x_j)$ whenever $i\geq j$. Any such sequence $X$
satisfies the conditions needed to serve as the generating sequence of
a {\em weighted pc-presentation} of $G$. The 
relations, $R$, in such a presentation all have the form
\begin{equation}
\label{eq:define-pc-relations}
\begin{array}{ccll}
x_i^p & = & \prod_{k=i+1}^n x_k^{b(i,k)} & 0\leq b(i,k)<p,~1\leq i\leq n \\
x_j^{x_i} & = & x_j\prod_{k=j+1}^nx_k^{b(i,j,k)} & 0\leq b(i,j,k)<p,~1\leq i<j\leq n
\end{array}
\end{equation}
We write $\la X\mid R\ra$ to denote the $p$-group defined by such a presentation.
We adopt the usual convention that an omitted relation $x_i^p$ implies
that $x_i^p=1$, and an omitted relation $x_j^{x_i}$ implies that $x_i$ and $x_j$
commute. We will often find it convenient to write a conjugate relation $x_j^{x_i}=x_jw$
as a commutator relation $[x_j,x_i]=w$.

\begin{remark}
{\rm
In general,
one requires that $G=P_1>\ldots>P_n>1$ refines a related series
called the {\em exponent $p$-central series}~\cite[p. 355]{handbook}. 
For the families of $p$-groups we consider here, however, the two series coincide.
}
\end{remark}


A critical feature of a pc-presentation for a $p$-group is that elements of the group inherit  
a {\em normal form} $x_1^{a_1}x_2^{a_2}\ldots x_n^{a_n}$, where $0\leq a_i<p$.
Given $g\in G$ as a word in $x_1,\ldots,x_n$, a normal form may be obtained by
repeatedly applying the relations in Equation~(\ref{eq:define-pc-relations}) in a process
known as ``collection".
If each element of $G$ has a unique normal form, the pc-presentation is said to be 
{\em consistent}. Clearly if $G$ has a consistent pc-presentation on $X=[x_1,\ldots,x_n]$,
then $|G|=p^n$. 

We conclude this section with a useful test for
consistency. We state it just for 2-groups -- since this is all we need --
and refer the reader to~\cite[Theorem 9.22]{handbook} for the more general version.

\begin{proposition}
\label{prop:consistency}
A weighted pc-presentation of a $d$-generator 2-group of class $c$ 
on $[x_1,\ldots,x_n]$ is consistent if the following pairs of words in the 
generators have the same normal form:
$$
\begin{array}{ll}
(x_kx_j)x_i~\mbox{and}~x_k(x_jx_i) & 1\leq i<j<k\leq n\;\mbox{and}\;i\leq d,\;w(x_i)+w(x_j)+w(x_k)\leq c; \\
(x_jx_j)x_i~\mbox{and}~x_j(x_jx_i) & 1\leq i<j\leq n\;\mbox{and}\;i\leq d,\;w(x_i)+w(x_j)<c; \\
(x_jx_i)x_i~\mbox{and}~x_j(x_ix_i) & 1\leq i<j\leq n,\;w(x_i)+w(x_j)<c; \\
(x_ix_i)x_i~\mbox{and}~x_i(x_ix_i) & 1\leq i\leq n,\;2w(a_i)<c;
\end{array}
$$
(The products in parentheses are collected first in each word.)
\end{proposition}

\section{The families \boldmath$\mathcal{H}^{\epsilon}$}
\label{four-families}
In this section we introduce four infinite families of 4-generator 2-groups of fixed coclass 4.
In the next section we will show that each family consists of groups that have
a class-preserving outer automorphism, thus proving Theorem~\ref{thm:main}.

We will define the groups in each family by giving consistent pc-presentations.
It is convenient to denote the ordered list of pc-generators of the 
$n^{{\rm th}}$ group in each
family by $X_n=\{x_1,x_2,x_3,x_4,z,y_1,\ldots ,y_n\}$,
with the group minimally generated by $\{x_1,x_2,x_3,x_4\}$. The commutator
relations for each family are identical, namely
\begin{equation}
\label{eq:comms}
C_n=\left\{
\begin{array}{ll}
 ~\hspace*{-4mm} & [x_2,x_1]=[x_3,x_2]=[x_4,x_1]=z,\;[x_3,x_1]=y_1,  
 \\ & [x_1,y_{i}]=[x_3,y_i]=y_{i+1}\;(i=1,\ldots,n-1) 
\end{array} 
\right\}
\end{equation}
For each
$\epsilon=(\epsilon_1,\epsilon_2,\epsilon_3,\epsilon_4)\in\{0,1\}^4$,
define
\begin{equation}
\label{eq:relns}
P_n^{\,\epsilon}=\left\{
\begin{array}{l}
x_j^2=z^{\epsilon_j}\;(j=1,\ldots,4),~z^2=1,~y_n^2=1\\ 
y_i^2=y_{i+1}y_{i+2}\;(i=1,\ldots,n-2),\;y_{n-1}^2=y_n
\end{array}
\right\}
\end{equation}
Let $R_n^{\,\epsilon}=C_n\cup P_n^{\,\epsilon}$,
define $H_n^{\epsilon}=\la X_n\mid R_n^{\,\epsilon}\ra$, and put
$\mathcal{H}^{\epsilon}=\{H_n^{\epsilon}\}_{n=1}^{\infty}$. Note that the pc-presentations
for the $n^{{\rm th}}$ group in each family differ only in the power relations of the
generators $x_i$.

\begin{proposition}
\label{prop:basic}
Let $n$ be a positive integer, and $\epsilon\in\{0,1\}^4$. Then $H_n^{\epsilon}=\la X_n\mid 
R_n^{\epsilon}\ra$ has order $2^{n+5}$ and class $n+1$ (hence coclass 4).
\end{proposition}

\noindent{\bf Proof.}
To confirm the order of $H_n^{\epsilon}$,
it suffices to check that their defining pc-presentations
are consistent, for which we use Proposition \ref{prop:consistency}.
Although there are $O(n^3)$ computations involved in
that test, the lion's share of these may be treated uniformly for the groups $H_n^{\epsilon}$.
The following table lists all of the 
triples that must be checked, together with their normal forms.
Triples involving $z$ are omitted, as this generator is central, as are
triples involving two or more $y_s$ generators, since $\la y_s\colon
s=1,\ldots,n\ra$ is abelian.
\begin{center}
\begin{tabular}{||c|c|c||}
\hline\hline
Triple $(a,b,c)$ & Conditions & Normal form of $a(bc)$ and $(ab)c$ \\
\hline\hline
$(x_3,x_2,x_1)$ &  & $x_1x_2x_3y_1$ \\
$(x_4,x_2,x_1)$ &  & $x_1x_2x_4$ \\
$(x_4,x_3,x_1)$ &  & $x_1x_3x_4zy_1$ \\
$(x_4,x_3,x_2)$ &  & $x_2x_3x_4z$ \\
\hline
$(y_s,x_2,x_1)$ & $s\leq n-2$ & $x_1x_2zy_sy_{s+1}$ \\
$(y_s,x_3,x_1)$ & $s\leq n-2$ & $x_1x_3y_1y_s$ \\
$(y_s,x_4,x_1)$ & $s\leq n-2$ & $x_1x_4zy_sy_{s+1}$ \\
$(y_s,x_3,x_2)$ & $s\leq n-2$ & $x_2x_3zy_sy_{s+1}$ \\
$(y_s,x_4,x_2)$ & $s\leq n-2$ & $x_2x_4y_s$ \\
$(y_s,x_4,x_3)$ & $s\leq n-2$ & $x_3x_4y_sy_{s+1}$ \\
\hline
$(x_j,x_j,x_i)$ & $1\leq i<j\leq 4$ & $x_iz^{e_j}$ \\
$(y_s,y_s,x_i)$ & $s\leq n-2$, $i=1,3$ & $x_iy_{s+1}$ \\
$(x_j,x_i,x_i)$ & $1\leq i<j\leq 4$ & $x_jz^{e_i}$ \\
$(y_s,x_i,x_i)$ & $s\leq n-2$, $i\leq 4$ & $z^{e_i}y_s$ \\
\hline
$(x_i,x_i,x_i)$ & $i\leq 4$ & $x_iz^{e_i}$ \\ \hline\hline
\end{tabular}
\end{center}

\noindent 
Routine calculations using the pc-relations are all that is needed to verify the normal forms listed in the table. 
It remains to compute the lower central series of $H_{n}^{\epsilon}$:
\[
\begin{array}{rcl}
\gamma_1(H_{n}^{\epsilon}) &=& H_{n}^{\epsilon} \\

\gamma_2 (H_{n}^{\epsilon}) &=&\langle z, y_i \colon 1\leq i\leq n \rangle\\

\gamma_j(H_n^{\epsilon}) & = & \langle y_{i} \colon j-1\leq i\leq n \rangle~~~\mbox{for}\;j=3,\ldots,n+1 \\

\gamma_{n+2}(H_{n}^{\epsilon}) &=& 1.
\end{array}
\]
This shows that $H_n^{\epsilon}$ has class $n+1$, as stated. ~$\Box$
\vspace*{3mm}

 Proposition~\ref{prop:basic} suggests that there are 16 families $\mathcal{H}^{\epsilon}$, but
 the following result shows that there is some duplication.
 
\begin{proposition}
For each positive integer $n$, there are four isomorphism classes among
the groups $\{H_n^{\epsilon}\colon\epsilon\in\{0,1\}^4\}$.
\end{proposition}

\noindent {\bf Proof.} Each group $H=H_n^{\epsilon}$ determines
a quadratic map $\qq=\qq^{\epsilon}$ (independent of $n$) as follows. 
Let $V$ denote
the largest elementary abelian quotient of $H$, namely
$V=H/A\cong(\Z/2)^4$, where $A=\la z,y_1,\ldots,y_n\ra$. Let $W$ denote the
largest elementary abelian quotient of $A$, namely $W=A/B\cong(\Z/2)^2$,
where $B=\la y_2,\ldots,y_n\ra$.
Define maps 
$\qq\colon V\to W$ and $\bb\colon V\times V\to W$,
where $\qq(xA)=x^2B$ and $\bb(xA,yA)=[x,y]B$ for all $x,y\in H$.
Using additive notation in $V$ and $W$, one easily checks that
$$
\bb(u,v)=\qq(u+v)+\qq(u)+\qq(v)~\mbox{for all}\;u,v\in V,
$$
so $\bb$ is the symmetric bilinear map associated to $\qq$ in
the familiar sense.

If $H_n^{\epsilon}$ and $H_n^{\delta}$ are isomorphic groups,
and $\alpha\colon H_n^{\epsilon}\to H_n^{\delta}$ is any isomorphism,
then $\alpha$ induces isomorphisms $\beta\colon V^{\epsilon}\to V^{\delta}$
and $\gamma\colon W^{\epsilon}\to W^{\delta}$ such that
$\qq^{\delta}(v\beta)=\qq^{\epsilon}(v)\gamma$ for all $v\in V^{\epsilon}$.
Thus $\alpha$ induces a {\em pseudo-isometry} between $\qq^{\epsilon}$
and $\qq^{\delta}$.

Given two such quadratic maps $\qq^{\epsilon}$ and $\qq^{\delta}$ one can easily
test for pseudo-isometry as follows. Fixing a basis $\{v_i\}$ for $V$, 
represent a quadratic map $\qq$ as  a $4\times 4$ matrix 
$\QQ=[[q_{ij}]]$ with entries in $W$, where $q_{ii}=\qq(v_i,v_i)$,
$q_{ij}=\bb(v_i,v_j)$ if $i<j$, and $q_{ij}=0$ of $i>j$. Then
$\qq(v)=v\QQ v^{{\rm tr}}$ for all $v\in V$. Using the basis
$\{x_iA\}$ for $V$, and identifying $A/B$ with the additive group of the ring
$(\Z/2)[t]/(t^2)$, the matrix representing $\qq$ is
$$
\QQ=\left[ \begin{array}{cccc} \epsilon_1 & 1 & t & 1 \\ 0 & \epsilon_2 & 1 & 0 \\
0 & 0 & \epsilon_3 & 0 \\ 0 & 0 & 0 & \epsilon_4 \end{array} \right],
$$
and the matrix representing the associated bilinear map $\bb$ is
$\BB=\QQ+\QQ^{{\rm tr}}$.

Let $\QQ^{\epsilon}$ and $\QQ^{\delta}$ be matrices representing the
quadratic maps associated to the groups $H^{\epsilon}$ and $H^{\delta}$
for $\epsilon,\delta\in\{0,1\}^4$. 
If $g\in{\rm GL}(4,2)$ represents an
isomorphism $H^{\epsilon}/A^{\epsilon}\to H^{\delta}/A^{\delta}$
induced by an isomorphism $H^{\epsilon}\to H^{\delta}$,
then the induced isomorphism $A^{\epsilon}/B^{\epsilon}\to A^{\delta}/B^{\delta}$
is uniquely determined by $g$, and its matrix $h\in{\rm GL}(2,2)$ is easily computed.
Extend $h$ entry-wise to a map 
$\Bbb{M}_4(W^{\epsilon})\to\Bbb{M}_4(W^{\delta})$, and denote
the image of $X\in\Bbb{M}_4(W^{\epsilon})$ by $X^{h}$.
Then $\qq^{\epsilon}$ and $\qq^{\delta}$ are pseudo-isometric if
and only if there exists $g\in{\rm GL}(4,2)$ such that
$$
g\BB^{\delta}g^{{\rm tr}}=(\BB^{\epsilon})^{h}~~\mbox{and}~~
v_i(g\QQ^{\delta}g^{{\rm tr}})v_i^{{\rm tr}}=v_i(\QQ^{\epsilon})^{h}v_i^{{\rm tr}},
$$
as $v_i$ runs over a basis for $(\Z/2)^4$.

Thus, the determination of the pseudo-isometry classes of the quadratic maps
associated to the families $\mathcal{H}^{\epsilon}$ is an elementary matrix
calculation in ${\rm GL}(4,2)$, which is easily carried out in {\sc magma}.
Those classes are represented by 
$$
\QQ^{\epsilon}~~\mbox{for}~~\epsilon\in\{(0,0,0,0),(0,1,0,0),(0,1,1,0),(0,0,0,1)\}.
$$

Finally, it is not difficult to verify that any pseudo-isometry $\QQ^{\epsilon}\to\QQ^{\delta}$
lifts to an isomorphism $H^{\epsilon}\to H^{\delta}$. Thus, for each $n$, there are precisely
four isomorphism classes of group $H_n^{\epsilon}$, as claimed.
 ~$\Box$
 \bigskip
 
 \noindent {\bf Remark.}~ The individual groups in each family can be constructed directly in
 {\sc magma} using their given pc-presentations. Alternatively, they can be constructed sequentially
 using the $p$-group generation algorithm. For example, 
 \smallskip
 
 \begin{tabular}{ll}
  $H2\;:=\;\mathtt{SmallGroup}\,(64,215)$; & \\
  $H3\;:=\;\mathtt{Descendants}\,(H2)[1]$;  & /* there are 19 descendants; $H3$ is the first */ \\
  $H4\;:=\;\mathtt{Descendants}\,(H3)[1]$;  & /* again, 19 descendants */ \\
 \end{tabular}
\smallskip

\noindent and so forth, constructs the members in the family $\mathcal{H}^{(0,0,0,0)}$.
Each group in the family 19 descendants, but only the first descendant in the returned 
list is ``capable" in the sense that it has further descendants. 
It was in precisely this way that the four families of groups were first discovered.
The remaining three families can be constructed by starting
instead with $\mathtt{SmallGroup}\;(64,n)$ for $n\in\{216,217,218\}$.

\section{Proof of Theorem~\ref{thm:main}}
\label{proof}
In this section we complete the proof of Theorem~\ref{thm:main} by exhibiting
a nearly inner automorphism of each group $H_n^{\epsilon}$ that is not inner.

Fix $n\geq 1,\; \epsilon\in\{0,1\}^4$, and put $H=H_n^{\epsilon}$. 
Define $\theta\colon H\to H$,
sending $x_4\mapsto x_4z$, and fixing all of the other generators in $X_n$.
One easily verifies (replacing $x_4$ by $x_4z$ in each pc-relation involving $x_4$ 
and evaluating) that $\theta\in\aut(H)$.

\begin{lemma}
\label{lem:not-inner}
The automorphism $\theta$ is not inner.
\end{lemma}

\noindent {\bf Proof.}
If $\theta$ is an inner automorphism, then there exists $h\in H$ commuting with
$x_1$ and $x_3$, but not with $x_4$. If $h=\prod_{i=1}^4x_i^{a_i}\cdot z^b\cdot\prod_{j=1}^ny_j^{c_j}$,
where all exponents are 0 or 1, then using the defining commutator relations of $H$ we see that
$$
hx_1=x_1h\cdot\left(z^{a_2+a_4}y_1^{a_3}\prod_{j=2}^ny_j^{c_{j-1}}\right).
$$
Hence $h\in C_H(x_1)$ if and only if
$a_2=a_4$ and $0=a_3=c_1=\ldots=c_{n-1}$. Also,
$$
\begin{array}{rcl}
x_3h & = & x_1^{a_1}x_2^{a_2}x_3^{1+a_3}x_4^{a_4}z^{a_2+b}y_1^{a_1+c_1}\prod_{j=2}^ny_j^{c_{j}},~\mbox{while} \\
hx_3 & = & x_1^{a_1}x_2^{a_2}x_3^{1+a_3}x_4^{a_4}z^by_1^{c_1}
\prod_{j=2}^ny_j^{c_{j}} \prod_{j=2}^ny_j^{c_{j-1}},
\end{array}
$$
so that $h\in C_H(x_3)$ if and only if $0=a_1=a_2=c_1=\ldots=c_{n-1}$. It follows that
$C_H(x_1)\cap C_H(x_3)=\la z,y_n\ra=Z(H)$. Thus $\theta$ is not inner. $\Box$
\vspace*{3mm}

The next lemma completes the proof of Theorem~\ref{thm:main}.

\begin{lemma}
\label{lem:nearly-inner}
The automorphism $\theta$ is nearly inner.
\end{lemma}

\noindent {\bf Proof.} 
We must show that, for each $h\in H$, there exists $t=t(h)\in H$ such that
$h^t=h\theta$. Fix $h\in H$, and write $h=\prod_{i=1}^4x_i^{a_i}\cdot z^b\cdot\prod_{j=1}^ny_j^{c_j}$,
as in the proof of the previous lemma. If $a_4=0$, then $h\theta=h$ and we may choose
$t(h)=1$. Thus, we may assume that $a_4=1$, whence $h\theta=hz$.
We claim that either $h^{x_2}=hz$ or $h^{x_1x_3}=hz$.

It is clear from the pc-relations that $x_2$ commutes with every $y_j$.
This is true also of $x_1x_3$. For, if $j<n-1$, then $y_j^{x_1x_3}=(y_jy_{j+1})^{x_3}=
y_jy_{j+1}^2y_{j+2}$. Using the relations (and a finite induction) one sees that
$y_{j+1}^2y_{j+2}=y_{n-1}^2y_n=y_n^2=1$. It is easy to see that $y_{n-1}^{x_1x_3}=y_{n-1}$
and that $y_n^{x_1x_3}=y_n$.

Next, observe that $x_2$ commutes with $x_4$, while $x_4^{x_1x_3}=(x_4z)^{x_3}=x_4z$.
Thus, it suffices to show, for each $(a_1,a_2,a_3)\in\{0,1\}^3$, if $h=x_1^{a_1}x_2^{a_2}x_3^{a_3}$,
then either $h^{x_2}=hz$, or $h^{x_1x_3}=h$. First,
$$
h^{x_2}=(x_1^{a_1}x_2^{a_2}x_3^{a_3})^{x_2}=x_1^{a_1}x_2^{a_2}x_3^{a_3}z^{a_1+a_3}=hz^{a_1+a_3}.
$$
Hence, if $a_1\neq a_3$, then $h^{x_2}=hz$, as required. 
It remains to show that $x_1x_3$ commutes with $h$ whenever $a_1=a_3$.
If $a_1=a_3=0$, then either $h=1$ or $h=x_2$;  clearly $x_1x_3$ commutes with 1, 
and $x_2^{x_1x_3}=x_2z^2=x_2$. Finally, if $a_1=a_3=1$, then either $h=x_1x_3$
or $h=x_1x_2x_3$; clearly $x_1x_3$ commutes with itself, and
$$
\begin{array}{rcl}
(x_1x_2x_3)^{x_1x_3} & = & (x_1(x_2z)(x_3y_1))^{x_3} \\
& = & (x_1y_1^{-1})(x_2z)zx_3(y_1y_2) \\
& = & x_1x_2y_1^{-1}x_3y_1y_2 \\
& = & x_1x_2x_3y_2^{-1}y_1^{-1}y_1y_2~=~x_1x_2x_3.
\end{array} 
$$
This completes the proof of the lemma. $\Box$

\section{Experimental data}
\label{data}

Using the 
{\sc magma} system we computed, for $n<512$, all groups $G$ of order $n$ such that 
${\rm Out_c}(G)\neq 1$. The results are summarized in the table below.
\begin{enumerate}
\item The first column of the table lists those orders for which there exists at least one 
group with a class-preserving outer automorphism. 
\item For each $n$, $d_n$ is the number of pairwise non-isomorphic groups of order $n$. 
\item If $C_n$ is the set of (isomorphism classes of) groups of order $n$
that possess a class-preserving outer automorphism, then $c_n$ is the
cardinality of $C_n$.
\item The fourth column records the set 
$O_n=\{ |{\rm Out_c}(G) |\colon G\in C_n\}$.
\end{enumerate}

\begin{tabular}[b]{ || b{10mm} | m{20mm} | m{20mm} | m{40mm}||}
\hline\hline
$n$ & $d_n$  & $c_n$ & $O_n$ \\ \hline
32 & 51 & 2 & \{ 2 \} \\
64 & 267 & 40 & \{ 2, 4, 16 \} \\
96 & 231 & 8 & \{ 2 \} \\
128 & 2,328 & 767 & \{ 2, 4, 8, 16, 64 \} \\
160 & 238 & 8 & \{ 2 \} \\
192 & 1,543 & 233 & \{ 2, 4, 16 \} \\
200 & 52 & 1 & \{ 2 \} \\
224 & 197 & 8 & \{ 2 \} \\
243 & 67 & 8 & \{ 3 \} \\
256 & 56,092 & 34,112 & \{ 2, 4, 8, 16, 32, 64 \} \\
288 & 1,045 & 28 & \{ 2 \} \\
300 & 49 & 1 & \{ 2 \} \\
320 & 1,640 & 243 & \{ 2, 4, 16 \} \\
352 & 195 & 8 & \{ 2 \} \\
384 &         &     &           \\
400 & 221 & 5 & \{ 2 \} \\
416 & 235 & 8 & \{ 2 \} \\
448 & 1,396 & 231 & \{ 2, 4, 16 \} \\
480 & 1,213 & 32 & \{ 2 \} \\
486 & 261 & 12 & \{ 3 \} \\ \hline\hline
\end{tabular}
\\
\\

The code used to conduct our experiment is available from either author. Briefly, the
method we used to compute ${\rm Out}_c(G)$ for a given group $G$ proceeds as follows.
\begin{enumerate}
\item Compute ${\rm Aut}(G)$. This is usually the most expensive step, particularly 
for certain soluble groups $G$. Indeed the default {\sc magma} function was incapable
of handling all of the groups of order 384 and we are indebted to David Howden
for supplying new code that enabled us to complete our search.
\item Compute $\Omega$, the set conjugacy classes of $G$, together
with the natural action $\rho\colon{\rm Aut}(G)\to {\rm Sym}(\Omega)$.
This step uses standard {\sc magma} functions.
\item Compute ${\rm ker}\rho$ and ${\rm Inn}(G)$ from which we construct
${\rm Out}_c(G)={\rm ker}\rho/{\rm Inn}(G)$. This step is fairly standard as well,
although we first convert from the usual representation of ${\rm Aut}(G)$ as
mappings to the standard (faithful) permutation representation of this group
on the elements of $G$.
This allows us to use the very efficient permutation group machinery in 
{\sc magma} to compute $\ker\rho$ quickly.
\end{enumerate}

There are over 12.5 million groups of order 512, which is why we chose to stop
at order 511. It would not be difficult, however, given sufficient computing power, to
extend the search to order 512 and beyond.
\medskip

We conclude this section, and the paper, by reporting on the type of groups
constructed by Heineken, namely those groups $G$ for which
${\rm Aut}(G)={\rm Aut}_c(G)$. 
There are precisely two such groups having order $<512$,
one having order 128 and the other having order 486.
\bigskip

\noindent {\bf Acknowledgments.} The authors would like to thank R. Quinlan
for bringing this problem to our attention, and E.A. O'Brien  for suggesting
the exhaustive search that led to the results in this paper.


\thispagestyle{plain}

\vspace*{10mm}

{\obeylines
Peter A. Brooksbank
Department of Mathematics
Bucknell University
Lewisburg, PA 17837
email: pbrooksb@bucknell.edu
\vspace*{7mm}

Matthew S. Mizuhara
Department of Mathematics
Bucknell University
Lewisburg, PA 17837
email: msm030@bucknell.edu
}

\end{document}